\documentclass[11pt,reqno]{amsart}
\usepackage{amsfonts,amsmath,amssymb,amsthm}
\usepackage[a4paper,text={6.1in,8.6in},centering]{geometry}
\usepackage{enumitem}
\usepackage[colorlinks=true,linkcolor=blue,citecolor=blue]{hyperref}
\usepackage{comment}
\usepackage{hyperref}
\usepackage{xcolor}

\theoremstyle{plain}
\newtheorem{thm}{Theorem}[section]

\newtheorem{lem}[thm]{Lemma}
\newtheorem{prop}[thm]{Proposition}
\theoremstyle{definition}
\newtheorem{assump}[thm]{Assumption}
\theoremstyle{remark}
\newtheorem*{remark}{Remark}
\numberwithin{equation}{section}

\newcommand{\LL}{\mathcal{L}}
\newcommand{\Sone}{\mathbb{S}^1}
\newcommand{\D}{\mathbb{D}}
\newcommand{\R}{\mathbb{R}}
\newcommand{\Z}{\mathbb{Z}}
\newcommand{\PP}{\mathbb{P}}
\newcommand{\EE}{\mathbb{E}}
\newcommand{\tflat}{\mathrm{tr}^{\flat}}
\newcommand{\Id}{\mathrm{Id}}

\title[Random zeta functions under a spectral gap]{Random zeta functions under a spectral gap}

\author[Y.\ Nakano]{Yushi Nakano}
\address{Department of Mathematics, Faculty of Science, Hokkaido University, Sapporo 060-0810, Japan}
\email{yushi.nakano@math.sci.hokudai.ac.jp}

\author[S.\ Suzuki]{Shintaro Suzuki}
\address{Department of Mathematics, Tokyo Gakugei University, Koganei, Tokyo 184-8501, Japan}
\email{shin05@u-gakugei.ac.jp}

\subjclass[2020]{Primary 37C30, 37H15; Secondary 37D20, 37A50}
\keywords{Random dynamical systems, dynamical zeta function, transfer operator, natural boundary, spectral gap}
\date{}

\begin{document}
\begin{abstract}
We study the meromorphic continuation of sample-wise zeta functions for i.i.d. random dynamical systems generated by $C^{r+1}$ local diffeomorphisms on a compact Riemannian manifold with $r\geq1$. Under general spectral gap assumptions on transfer operator cocycles and suitable conditions on their flat traces, we establish two results. First, we prove that the sample-wise zeta function associated with the geometric potential admits a meromorphic continuation. Second, under further assumptions, we show that the sample-wise zeta function has a natural boundary, providing an analogue of Buzzi's
result for random products of $2\times 2$ matrices. For uniformly expanding families of $C^{r+1}$ maps, we verify the spectral gap and flat-trace assumptions used in these results.
\end{abstract}

\maketitle

\section{Introduction}\label{sec:intro}

\subsection{Background}

Given a $C^{r+1}$ expanding map $T$ on a compact Riemannian manifold $M$ and a positive weight $g\in C^r(M)$ with $r\geq1$, the dynamical zeta function
\[
\zeta(z):=\exp\Bigl(\sum_{n\geq 1}\frac{z^n}{n}\sum_{T^n x=x}g^{(n)}(x)\Bigr)
\]
converges in a certain disk and admits a meromorphic continuation beyond its circle of convergence. This is often proved using dynamical determinants via the transfer operator method; see, for example, Ruelle \cite{Rue89} and Baladi~\cite{Bal00}.

For independent and identically distributed (i.i.d.) random dynamical systems, Ruelle~\cite{Rue89} showed that, in the expanding case, the corresponding sample-mean (annealed) zeta function also converges in a certain disk and admits a meromorphic continuation beyond its circle of convergence under suitable conditions.
Baladi~\cite{Bal98} raised the question of whether the same holds for sample-wise (quenched) zeta functions almost surely, as these functions can be expected to encode information about transfer operator cocycles.
Buzzi~\cite{Buz02} showed that this is not true in general, even for random products of $2\times 2$ matrices: almost surely, the sample-wise zeta function has the unit circle as a natural boundary, in the sense that it admits no meromorphic continuation across any point of the circle.

In this paper, we consider i.i.d. random dynamical systems generated by $C^{r+1}$ local diffeomorphisms on a compact Riemannian manifold and 
establish two results concerning the meromorphic extension of sample-wise zeta functions. The first result contrasts with Buzzi's result, while the second provides an analogue of Buzzi's result. Our framework is formulated in terms of general spectral gap assumptions on transfer operator cocycles, together with mild conditions on their flat traces. 
In Section~\ref{sec:example}, we verify the spectral gap and flat-trace assumptions for uniformly expanding families of maps with uniformly positive weights. In this setting, the first result applies to the geometric potential, while the second requires only a non-degeneracy condition.
The precise notation used in the following statements is introduced in Section~\ref{sec:setting}.
\begin{thm}\label{thm:mero}
Assume \hyperref[ass:SG1]{(SG1)} from Section~\ref{sec:setting}.
Then the random zeta function $\zeta_\omega(z)$ for the geometric potential extends $\PP$-almost surely meromorphically to
\[
|z|<\min\{\rho^{-1},\alpha\},
\]
where $\rho\in(0,1)$ and $\alpha>1$ are the constants in \hyperref[ass:SG1]{(SG1)}.
\end{thm}
For the geometric potential, the leading Lyapunov exponent
$\EE[\log\lambda]$ of the transfer operator cocycle vanishes.
Thus Theorem~\ref{thm:mero} gives meromorphic continuation across
the circle $|z|=\exp(-\EE[\log\lambda])$, which plays a very
different role in our second result.
\begin{thm}\label{thm:boundary}
Assume \hyperref[ass:SG2]{(SG2)} from Section~\ref{sec:setting}.
Then the circle $$|z|=\exp(-\EE[\log\lambda])$$ is $\PP$-almost surely a natural boundary for meromorphic continuation of the random zeta function $\zeta_\omega(z)$, where $\lambda$ is the random positive number in \hyperref[ass:SG2]{(SG2)}.
\end{thm}
We make some comments on assumptions \hyperref[ass:SG1]{(SG1)} and \hyperref[ass:SG2]{(SG2)} in Theorems \ref{thm:mero} and \ref{thm:boundary}, respectively. Both assumptions include \hyperref[ass:SG]{(SG)}, which consists of a spectral gap property for transfer operator cocycles and the condition that the difference between the weighted fixed-point sum appearing in the definition of sample-wise zeta functions and the corresponding flat trace is exponentially small relative to the leading Lyapunov eigenvalue $\lambda_\omega$.
The latter condition reduces the problem of meromorphic continuation of sample-wise zeta functions to that of dynamical determinants, allowing us to use the theory of flat traces.
In \hyperref[ass:SG1]{(SG1)}, we additionally assume that the weight is given by the the geometric potential.
This implies that all leading Lyapunov eigenvalues are equal to $1$ and that the sample-wise reference measures appearing in the spectral gap property are the normalized Lebesgue measure.
Consequently, the corresponding flat traces converge to $1$ exponentially fast, yielding Theorem \ref{thm:mero}.
In contrast, \hyperref[ass:SG2]{(SG2)} ensures that changing only the last map and its weight in a random composition
produces, with uniformly positive probability, a change in its flat trace bounded below by a
positive multiple of the leading Lyapunov eigenvalues, even for arbitrarily long compositions.

This paper is organized as follows. Section~\ref{sec:mero} establishes the meromorphic continuation result for the geometric potential. Section~\ref{sec:boundary} presents the argument establishing a natural boundary under some mild assumptions. Section~\ref{sec:example} verifies the spectral gap and flat-trace assumptions for uniformly expanding maps.

\section{Setting and assumptions}\label{sec:setting}

\begin{subsection}{Basic notions}\label{subsec:basic}
Throughout this paper, we consider a probability space $(M,\mathcal {M}, m)$, where $M$ is a compact boundaryless Riemannian manifold of dimension $d$, $\mathcal{M}$ is the Borel $\sigma$-algebra derived from the topology of $M$ and $m$ is the normalized Riemannian volume measure.
Let $(\Omega, \mathcal{F}, \mathbb P)=(\tilde\Omega, \tilde{\mathcal F}, \tilde{\mathbb P} )^{\mathbb{Z}}$ be the two-sided product probability space of some Lebesgue space $(\tilde\Omega, \tilde{\mathcal F}, \tilde{\mathbb P} )$ and let $\sigma:\Omega\to\Omega$ be the left shift, i.e., $\sigma(\omega_i)_{i\in\Z}=(\omega_{i+1})_{i\in\Z}$ for $(\omega_i)_{i\in\Z}\in\Omega$. Let $T:\tilde\Omega\to C^{r+1}(M,M)$ be a map such that $T_{\tilde{\omega}}:=T(\tilde\omega)$ for $\tilde{\mathbb{P}}$-a.e. $\tilde{\omega}\in\tilde{\Omega}$ is a ${C}^{r+1}$ local diffeomorphism on $M$ with $r\ge1$, uniformly bounded $ C^{r+1}$ norm, and uniformly bounded degree. Assume that the map $((\omega_i)_{i\in\Z},x)\in\Omega\times M\mapsto T_{\omega_0}(x)\in M$ is measurable. 
By a slight abuse of notation, we define $T_\omega=T_{\omega_0}$ for $\omega=(\omega_i)_{i\in\Z}\in\Omega$. We write 
\[T_\omega^n=T_{\omega_{n-1}}\circ\dots\circ T_{\omega_0}\]
for $n\geq0$ and $\omega=(\omega_i)_{i\in\Z}\in\Omega$, which is regarded as the identity map if $n=0$. Then the random dynamical system generated by $T$ is given by the skew-product map 
\[R(\omega,x)=(\sigma\omega, T_\omega(x))\]
for $(\omega,x)\in \Omega\times M$. 

Let $g:\tilde \Omega\to C^{r}(M,\mathbb R)$ be a function such that $g_{\tilde{\omega}}:=g(\tilde{\omega})$ is a positive $C^{r}$ function on $M$ for $\tilde\omega\in\tilde{\Omega}$ satisfying there is a positive constant $C_0>0$ with $||g_{\tilde\omega}||_{C^r}\leq C_0$ for any $\tilde\omega\in\tilde\Omega$ and assume that the map $((\omega_i)_{i\in\Z},x)\in\Omega\times M\mapsto g_{\omega_0}(x)\in\mathbb R$ is measurable. We define $g_\omega=g_{\omega_0}$ for $\omega=(\omega_i)_{i\in\Z}\in\Omega$ and set 
\[
 g_\omega^{(n)}=\prod_{j=0}^{n-1} g_{\omega_j}\circ T_\omega^j
\]
for $n\geq1$ and $((\omega_i)_{i\in\Z},x)\in\Omega\times M$. 

Let $g$ be a function as defined above. We define the weighted transfer operator by 
\[
(\LL_\omega u)(x)=\sum_{T_\omega y=x} g_\omega(y)u(y)
\]
for $u\in C^{\infty}(M,\mathbb R)$ and $x\in M$. Note that $g_\omega$ and $T_\omega$ for $\omega\in\Omega$ only depends on $0$-th coordinate of $\omega$. We set $\mathcal L_\omega=\mathcal L_{\omega_0}$ for $\omega=(\omega_i)_{i\in\Z}\in\Omega$ and write
\[
\LL_\omega^n:=\LL_{\sigma^{n-1}\omega}\cdots\LL_\omega=\mathcal L_{\omega_{n-1}}\cdots\mathcal L_{\omega_0}
\]
for the corresponding operator cocycle, where $\omega=(\omega_i)_{i\in\Z}$ and $n\geq1$. We define the weighted fixed-point sum by 
\[
Z_n(\omega):=\sum_{T_\omega^n x=x} g_\omega^{(n)}(x)
\]
for $n\geq1$ and $\omega\in\Omega$. We then define the sample-wise zeta function by
\[
\zeta_\omega(z)=\exp\Bigl(\sum_{n\ge1}\frac{z^n}{n}Z_n(\omega)\Bigr),
\]
regarded as a formal power series in $z$. 

According to \cite[Chapter 3]{Bal18}, we introduce a Banach space $\mathcal B$ on which transfer operators act and a flat trace for associated operators can be defined. 
For expanding maps, spaces of differentiable or Sobolev functions on $M$ are natural choices.
We assume that there are a Banach space $(\mathcal B, ||\cdot ||_{\mathcal B})$ in which $C^\infty(M,\R)$ is dense with respect to the norm $||\cdot||_\mathcal{B}$ and operators $I_\varepsilon: \mathcal{B}\to \mathcal{B}$ for $\varepsilon\in(0,\varepsilon_0)$, where $\varepsilon_0>0$ is a positive constant, satisfying the following conditions: 

\begin{itemize}

\item[(1)] For $\varepsilon\in(0,\varepsilon_0)$ there is a $C^\infty$ function $K_\varepsilon:M\times M\to [0,\infty)$ so that 
\[I_\varepsilon(\varphi)(\cdot)=\int_M K_\varepsilon(\cdot, y)\varphi(y)dy\]
for $\varepsilon\in(0,\varepsilon_0)$ and $\varphi\in C^\infty(M,\mathbb R)$. 

\item[(2)] The support of $K_\varepsilon$ tends to $\{(x,x)\in M\times M: x\in M\}$ as $\varepsilon\to0$. 

\item[(3)] $\sup_{y\in M}||K_\varepsilon(\cdot,y)||_{{C^t(M)}}<\infty$ for any $t\geq1$ and $\varepsilon\in(0,\varepsilon_0)$. 

\item[(4)] There is a positive constant $C>0$ such that $||I_{\varepsilon}(\varphi)||_{\mathcal B}\leq C||\varphi||_{\mathcal B}$
for any $\varepsilon\in(0,\varepsilon_0)$ and $\varphi\in\mathcal B$. 

\item[(5)] $\lim_{\varepsilon\to0}||I_{\varepsilon}(\varphi)-\varphi||_{\mathcal B}=0$
for any $\varphi\in \mathcal{B}$. 
\end{itemize}
If a bounded operator $S:\mathcal B\to \mathcal B$ satisfies 
\begin{equation}\label{T}
    S(C^t(M,\mathbb R))\subset C^0(M,\mathbb R) \text{ and } ||S\varphi||_{C^0}\leq C||\varphi||_{\mathcal{B}}
\end{equation}
for $t\geq1$ and $\varphi\in C^t(M,\mathbb R)$, where $C>0$ is a constant independent of $t\geq1$, then $S_\varepsilon(x,y):=S(K_\varepsilon(\cdot,y))(x)$ is in $C^0(M\times M,\mathbb R)$ for $x,y\in M$ and $\varepsilon\in(0,\varepsilon_0)$, which shows that the operator
\[SI_\varepsilon\varphi(x)=\int_M S_\varepsilon(x,y)\varphi(y)dy\]
for $\varphi\in\mathcal B$ and $\varepsilon\in(0,\varepsilon_0)$ is well-defined.
For a bounded operator $S: \mathcal{B}\to\mathcal{B}$ satisfying  (\ref{T}), we define 
\[\tflat(S)=\lim_{\varepsilon\to0}\tflat_\varepsilon(S)\]
if the limit exists, where $\tflat_\varepsilon(S)=\int_MS_\varepsilon(x,x)dx$. 
In the uniformly expanding setting of Section~\ref{sec:example}, the flat trace is given by the following formula (cf.~\cite[Proposition~3.14]{Bal18}):
\[
\tflat(\LL_\omega^n)
 =\sum_{T_\omega^n x=x}\frac{g_\omega^{(n)}(x)}{|\det(\Id-(D_xT_\omega^n)^{-1})|},
\qquad n\ge1.
\]

\end{subsection}

\begin{subsection}{The main assumptions}

In this section, we state the precise assumptions for the main results of the paper. Let $(T, g)$ be a pair of a map $T:\Omega\times M\to M$ generating an i.i.d. random dynamical system and a potential function $g:\Omega\times M\to \R$ as defined in subsection~\ref{subsec:basic}. 

\begin{assump}[SG]\label{ass:SG}
We say that $(T,g)$ satisfies \hyperref[ass:SG]{(SG)} if: 
\begin{itemize}
\item[(1)] There exist measurable families
$h_\omega \in \mathcal B$, $\nu_\omega \in \mathcal B^*$, and $\lambda_\omega>0$,
constants $C_1>0$ and $\rho\in(0,1)$, and a full-measure set
$\Omega_1\subset\Omega$ such that, for every $\omega\in\Omega_1$
and $n\ge1$, the flat trace of $\mathcal{L}_\omega^n$ is well defined and
\[
\mathcal{L}_\omega^n
=
\lambda_\omega^n h_{\sigma^n\omega}\otimes\nu_\omega
+
R_\omega^n,
\]
where $R_\omega^n:\mathcal B\to \mathcal B$ satisfies
\[
\left|\nu_\omega(h_{\sigma^n\omega})\right|\le C_1,
\qquad
\left\|R_\omega^n\right\|_B
\le C_1\rho^n\lambda_\omega^n,
\]
and
\[
\left|\operatorname{tr}^{\flat}(R_\omega^n)\right|
\le C_1\rho^n\lambda_\omega^n,
\]
where $\lambda_\omega^n=\prod_{i=0}^{n-1}\lambda_{\sigma^i\omega}$. 

\item[(2)] There are constants $C_{2}>0$, $\alpha>1$, and a full-measure set $\Omega_{\mathrm{2}}\subset\Omega$ such that, for every $\omega\in\Omega_{2}$ and every $n\ge1$,
\[
 |Z_n(\omega)-\tflat(\LL_\omega^n)|\le C_{2} \alpha^{-n}\lambda_\omega^n.
\]
\end{itemize}
\end{assump}

\begin{assump}[SG1]\label{ass:SG1}
We say that $(T,g)$ satisfies \hyperref[ass:SG1]{(SG1)} if it satisfies  \hyperref[ass:SG]{(SG)} and the potential $g$ is given by the geometric potential $|\det DT_\omega|^{-1}$. In particular,  
\[
\lambda_\omega\equiv1,
\qquad
\nu_\omega(u)=\int_M u\,dm,
\qquad
\int_M h_\omega\,dm=1.
\]
\end{assump}

\begin{assump}[SG2]\label{ass:SG2}
We say that $(T,g)$ satisfies \hyperref[ass:SG2]{(SG2)} if it satisfies  \hyperref[ass:SG]{(SG)} and the following hold: 
\begin{itemize}
\item [(1)] There are a constant $C_3>0$ and a full-measure set $\Omega_{3}\subset\Omega$ such that, for every $\omega\in\Omega_{3}$, every $n\ge1$, and for $\tilde{\mathbb{P}}$-almost every $\xi\in\widetilde\Omega$,
\[
\|(\LL_\xi-\LL_{\omega_n})R_\omega^n\|_{\mathcal B}
\le C_3\rho^n\lambda_\omega^n
\]
and 
\[
\left|\tflat\bigl((\LL_\xi-\LL_{\omega_n})R_\omega^n\bigr)\right|
\le C_3\rho^n\lambda_\omega^n.
\]

\item[(2)] There exist $\delta>0$, $p_0>0$, and a full-measure set $\Omega_{4}\subset\Omega$ such that, for every $\omega=(\omega_i)_{i\in\mathbb Z}\in\Omega_{4}$ and every $n\ge1$,
\[
\widetilde{\mathbb{P}}\Bigl\{\xi\in\widetilde\Omega:
\bigl|\nu_\omega\bigl((\LL_\xi-\LL_{\omega_n})h_{\sigma^n\omega}\bigr)\bigr|\ge 3\delta
\Bigr\}\ge p_0.
\]

\item[(3)] $\log \lambda$ is $\PP$-integrable.
\end{itemize}
\end{assump}

For $\omega\in\Omega$, $\xi\in \tilde{\Omega}$ and $i\in \Z$, we define 
\begin{equation}\label{Tn}
\mathcal{T}_n(\omega,\xi)_i
=\begin{cases}
    \omega_i, & (i\neq n), \\
    \xi, &(i=n).
\end{cases}
\end{equation}

\begin{prop}\label{prop:sd-consequences}
Assume that $(T,g)$ satisfies item $(1)$ in \hyperref[ass:SG]{(SG)}. Then 
\begin{equation}\label{ass:SGa}
\tflat(\LL_\omega^n)
=\lambda_\omega^n\nu_\omega(h_{\sigma^n\omega})+r_n(\omega),
\qquad
|r_n(\omega)|\le C_1\rho^n\lambda_\omega^n.
\tag{SGa}
\end{equation}
In addition, if item (1) in \hyperref[ass:SG2]{(SG2)} holds, then for $\mathbb P$-almost every $\omega\in\Omega_{1}\cap\Omega_{3}$, every $n\ge1$, and $\tilde{\mathbb{P}}$-almost every $\xi\in\tilde{\Omega}$,
\begin{equation}\label{ass:SGb}
\tflat(\LL_{\mathcal T_n(\omega,\xi)}^{n+1})-\tflat(\LL_\omega^{n+1})
=\lambda_\omega^n\nu_\omega\bigl((\LL_\xi-\LL_{\omega_n})h_{\sigma^n\omega}\bigr)+\varepsilon_n(\omega,\xi),
\tag{SGb}
\end{equation}
where
\[
|\varepsilon_n(\omega,\xi)|\le C_3\rho^n\lambda_\omega^n.
\]
\end{prop}

\begin{proof}
The flat-trace estimate for $R_\omega^n$ is part of item (1) in \hyperref[ass:SG]{(SG)}.
Since both $\mathcal L_\omega^n$ and $R_\omega^n$ have flat traces, their difference $\lambda_\omega^n h_{\sigma^n\omega}\otimes\nu_\omega$ satisfies \eqref{T}.
By the argument in the proof of \cite[Proposition~3.13]{Bal18}, the flat trace of a rank-one operator $u\otimes\ell$ satisfying \eqref{T} is $\ell(u)$.
Therefore
\[
\tflat(\lambda_\omega^n h_{\sigma^n\omega}\otimes\nu_\omega)
=\lambda_\omega^n\nu_\omega(h_{\sigma^n\omega}).
\]
Linearity of the flat trace gives \hyperref[ass:SGa]{(SGa)}.

For the formula \hyperref[ass:SGb]{(SGb)}, note that $\mathcal T_n(\omega,\xi)\in\Omega_1$ for almost every $(\omega,\xi)$, since $(\mathcal T_n)_*(\PP\otimes\widetilde{\PP})=\PP$ by independence.
Changing the $n$-th symbol gives
\[
\LL_{\mathcal T_n(\omega,\xi)}^{n+1}-\LL_\omega^{n+1}
=(\LL_\xi-\LL_{\omega_n})\LL_\omega^n.
\]
Insert the decomposition from item (1) in \hyperref[ass:SG]{(SG)}. Then the flat trace of $\lambda_\omega^n(\LL_\xi-\LL_{\omega_n})h_{\sigma^n\omega}\otimes\nu_\omega 
$
is
\[
\lambda_\omega^n
\nu_\omega\bigl((\LL_\xi-\LL_{\omega_n})h_{\sigma^n\omega}\bigr).
\]
The remaining term is $(\LL_\xi-\LL_{\omega_n})R_\omega^n$, whose flat trace is bounded by item (1) in \hyperref[ass:SG2]{(SG2)}.  This gives \hyperref[ass:SGb]{(SGb)}.
\end{proof}

A useful sufficient condition for item (2) in \hyperref[ass:SG2]{(SG2)} is the following.
\begin{prop}\label{prop:anti-concentration}
Suppose there exists $C_{\mathrm{ac}}>0$ and a full-measure set $\Omega_{\mathrm{ac}}\subset\Omega$ such that, for every $\omega\in\Omega_{\mathrm{ac}}$ and every $n\ge1$, the map
$
\xi\mapsto\nu_\omega(\LL_\xi h_{\sigma^n\omega})
$
is measurable and real-valued for all $\xi$, and its pushforward law under $\widetilde{\mathbb{P}}$ is absolutely continuous with respect to Lebesgue measure on $\R$, with density bounded by $C_{\mathrm{ac}}$.
Then item (2) in \hyperref[ass:SG2]{(SG2)} holds.
\end{prop}

\begin{proof}
Write $F_{\omega,n}(\xi):=\nu_\omega(\LL_\xi h_{\sigma^n\omega})$.
Let $\omega\in\Omega_{\mathrm{ac}}$ and $n\geq1$.
By the density bound,
\[
\widetilde{\mathbb{P}}\bigl(|F_{\omega,n}(\xi)-F_{\omega,n}(\omega_n)|<3\delta\bigr)
\le 6 C_{\mathrm{ac}} \delta.
\]
Choose $\delta>0$ so small that $6 C_{\mathrm{ac}} \delta<1$, and set $p_0:=1-6 C_{\mathrm{ac}}\delta>0$.
Then
\[
\widetilde{\mathbb{P}}\bigl(|F_{\omega,n}(\xi)-F_{\omega,n}(\omega_n)|\ge 3\delta\bigr)
\ge p_0,
\]
which is exactly item (2) in \hyperref[ass:SG2]{(SG2)}.
\end{proof}

\noindent\textit{Example.}
Let $T$ be a uniformly expanding $C^{r+1}$ map and $g$ a
strictly positive $C^r$ weight. Set
\[
T_\xi=T,\qquad g_\xi=g+\xi,
\]
where $\xi$ is uniformly distributed on $[0,A]$ for some $A>0$.
The map $\xi\mapsto\nu_\omega(\mathcal L_\xi h_{\sigma^n\omega})$
is affine with slope
$\nu_\omega(\sum_{Ty=\cdot}h_{\sigma^n\omega}(y))$.
By Lemma~\ref{lem:s5-spectral}, $h_\omega$ is uniformly
bounded above and away from zero, and $\nu_\omega$ is
positive with $\nu_\omega(h_\omega)=1$.
Since $T$ is surjective, the slope is therefore bounded below by
\[
\inf_M h_{\sigma^n\omega}\,\nu_\omega(1)
\ge \frac{\inf_M h_{\sigma^n\omega}}{\|h_\omega\|_{C^0}}
\ge c_0
\]
for some $c_0>0$ independent of $\omega$ and $n$.
Hence the law of $\xi \mapsto \nu_\omega(\mathcal L_\xi h_{\sigma^n\omega})$
has a density bounded by $(Ac_0)^{-1}$, and
Proposition~\ref{prop:anti-concentration} applies.

\begin{remark}
For the geometrical potential, the transfer operator preserves the integral with respect to the normalized Riemannian volume $m$. Hence
\[
\nu_\omega(\mathcal{L}_\xi h_{\sigma^n\omega})
=\int h_{\sigma^n\omega}\, dm
=
1.
\]
Therefore, item~(2) in \hyperref[ass:SG2]{(SG2)} fails identically.
\end{remark}

\end{subsection}

\section{Meromorphic continuation for the geometrical potential}\label{sec:mero}

Throughout the proofs, we write
\[
\tau_n(\omega):=\tflat(\LL_\omega^n),\qquad n\ge1
\]
for short.

\begin{proof}[Proof of Theorem~\ref{thm:mero}]
Under \hyperref[ass:SG1]{(SG1)},
\[
\tau_n(\omega)=\nu_\omega(h_{\sigma^n\omega})+r_n(\omega)=1+r_n(\omega),
\qquad |r_n(\omega)|\le C_1\rho^n
\]
for $n\geq1$. Therefore
\[
 d_{0,\omega}(z):=\exp\Bigl(-\sum_{n\ge1}\frac{z^n}{n}\tau_n(\omega)\Bigr)
\]
can be written as
\begin{equation}\label{d0}
 d_{0,\omega}(z)
 =\exp\Bigl(-\sum_{n\ge1}\frac{z^n}{n}(1+r_n(\omega))\Bigr)
 =(1-z)e^{G_\omega(z)},
\end{equation}
where
\[
G_\omega(z):=-\sum_{n\ge1}\frac{z^n}{n}r_n(\omega).
\]
Since $|r_n(\omega)|\le C_1\rho^n$, the function $G_\omega$ is holomorphic on $|z|<\rho^{-1}$.

Since \hyperref[ass:SG1]{(SG1)} yields the normalization $\lambda_\omega^n\equiv1$, the function 
\[
H_\omega(z):=\exp\Bigl(\sum_{n\ge1}\frac{z^n}{n}\,(Z_n(\omega)-\tau_n(\omega))\Bigr)
\]
is holomorphic on $|z|<\alpha$.
Since
\[
\zeta_\omega(z)=H_\omega(z)d_{0,\omega}(z)^{-1}
\]
as germs at the origin, the identity theorem gives the same identity on the common connected domain of continuation.
Since $d_{0,\omega}^{-1}$ is meromorphic on $|z|<\rho^{-1}$, the stated continuation follows.
\end{proof}

\section{Natural boundary via one-step perturbations}\label{sec:boundary}

We now prove Theorem~\ref{thm:boundary}.
Let $c=\exp(-\EE[\log\lambda])$.
Replacing $g_\xi$ by $cg_\xi$ replaces $\mathcal L_\xi$,
$\lambda_\omega$ and $R_\omega^n$ by $c\mathcal L_\xi$,
$c\lambda_\omega$ and $c^nR_\omega^n$, respectively,
with $h_\omega$ and $\nu_\omega$ unchanged.
The rescaled pair still satisfies \hyperref[ass:SG2]{(SG2)},
possibly with different constants, and the corresponding zeta function
is $\zeta_\omega(cz)$.
Since $\EE[\log(c\lambda)]=0$, it suffices to prove the theorem
under the normalization $\EE[\log\lambda]=0$.

Similar to Section~\ref{sec:mero} set $\tau_n(\omega)=\tflat(\LL_\omega^n)$ for $n\ge1$.
By Proposition~\ref{prop:sd-consequences},
\[
|\tau_n(\omega)|\le C_1(1+\rho^n)\lambda_\omega^n.
\]
Since $\EE[\log\lambda]=0$, Birkhoff's ergodic theorem gives
\[
\lim_{n\to\infty}\frac1n\log\lambda_\omega^n=0
\qquad\text{for $\PP$-almost every }\omega.
\]
Hence $\limsup_{n\to\infty}|\tau_n(\omega)|^{1/n}\le1$ for $\PP$-almost every $\omega$, and therefore
\[
Q_\omega(z):=\sum_{n\ge1}\tau_n(\omega)z^n
\]
is holomorphic on $\D$ for $\PP$-almost every sample.
The same argument, together with item (2) in Assumption~\ref{ass:SG}, gives $\limsup_{n\to\infty} |Z_n(\omega)|^{1/n}\le1$.
Hence both $\zeta_\omega$ and $d_{0,\omega}$ are well-defined and holomorphic on $\D$ for almost every sample.
Moreover $d_{0,\omega}$ is zero-free on $\D$ by its exponential definition.

\subsection{Two analytic lemmas}

\begin{lem}\label{lem:common-domain}
Assume that, on a set of positive probability, $d_{0,\omega}$ is holomorphic and zero-free on $\D$ and admits a meromorphic continuation across at least one point of $\Sone$.
Then there exist a simply connected domain $U\supset\D$, a conformal map $\phi\colon\D\to U$ with $\phi(0)=0$ and $L:=\phi'(0)>1$, a number $r\in(L^{-1},1)$, a compact set $K:=\phi(\{w:|w|\le r\})\Subset U$, a constant $M>0$, and a measurable set $\Omega_*\subset\Omega$ of positive probability such that, for every $\omega\in\Omega_*$,
\[
Q_\omega(z)=-z\frac{d'_{0,\omega}(z)}{d_{0,\omega}(z)}\]
extends holomorphically to $U$ and satisfies
\[
\sup_{z\in K}|Q_\omega(z)|\le M.
\]
\end{lem}

\begin{proof}
By the hypothesis and a countable base of disks with centres and radii in fixed countable dense sets, there exist a boundary point $\beta\in\Sone$, an $\varepsilon>0$, and a positive-probability set $E$ such that $d_{0,\omega}$ extends meromorphically to $\D\cup B(\beta,\varepsilon)$ for every $\omega\in E$.
The meromorphic continuation is unique on this connected set.
For each such $\omega$, the zeros and poles of this extension are isolated.
Hence there is a smaller ball $B(\gamma,s)\subset B(\beta,\varepsilon)$, with $\gamma$ in a fixed countable dense subset of $B(\beta,\varepsilon)\cap\Sone$ and $s>0$ rational, on which $d_{0,\omega}$ is holomorphic and non-vanishing.
The countable union of the corresponding subfamilies covers $E$, so one choice has positive probability.
After shrinking $E$ to this subfamily, we keep the notation $\gamma$ and write the radius as $\varepsilon_1>0$.
Then, for every $\omega\in E$, the function $d_{0,\omega}$ is holomorphic and non-vanishing on $B(\gamma,\varepsilon_1)$.
Hence
\[
Q_\omega(z)=-z\,d'_{0,\omega}(z)/d_{0,\omega}(z)
\]
is holomorphic on
\[
U:=\D\cup B(\gamma,\varepsilon_1).
\]
The intersection of the two disks is connected, so $U$ is simply connected.
Let $\phi\colon\D\to U$ be the Riemann map with $\phi(0)=0$ and $\phi'(0)>0$.
Since $U$ strictly contains $\D$, applying Schwarz's lemma to the inverse map restricted to $\D$ gives $L:=\phi'(0)>1$.
Choose any $r\in(L^{-1},1)$ and set $K:=\phi(\{w:|w|\le r\})\Subset U$.
For each $\omega\in E$, the quantity $\sup_K|Q_\omega|$ is finite.
Measurability follows by analytic continuation from $\mathbb D$.
Indeed, for each fixed $z\in U$, choose a finite chain of overlapping disks in $U$, independently of $\omega$, starting with a disk contained in $\mathbb D$ and ending with a disk containing $z$.
At each step, the Taylor coefficients of the continuation are given by Cauchy integrals of the preceding branch.
Hence $\omega\mapsto Q_\omega(z)$ is measurable for every fixed $z\in U$.
Taking the supremum over a countable dense subset of $K$ then shows that
$
\omega\longmapsto \sup_{z\in K}|Q_\omega(z)|
$
is measurable.
By the countable decomposition $E=\bigcup_{m\ge1}\{\omega\in E:\sup_K|Q_\omega|\le m\}$, some $M>0$ yields a subset $\Omega_*\subset E$ of positive probability with the required uniform bound.
\end{proof}

\begin{lem}\label{lem:coeff-decay}
In the setting of Lemma~\ref{lem:common-domain}, there exist constants $C>0$ and $\vartheta\in(0,1)$ such that the following holds.
If $\omega,\omega'\in\Omega_*$ satisfy
\[
\tau_j(\omega')=\tau_j(\omega)\qquad (1\le j\le n),
\]
then
\[
|\tau_{n+1}(\omega')-\tau_{n+1}(\omega)|\le C\vartheta^n.
\]
\end{lem}

\begin{proof}
Set
\[
F(z):=Q_{\omega'}(z)-Q_\omega(z)=\sum_{m\ge n+1}\bigl(\tau_m(\omega')-\tau_m(\omega)\bigr)z^m.
\]
Then $F$ is holomorphic on $U$ and satisfies $\sup_K|F|\le2M$.
Consider
\[
G(w):=F(\phi(w)),\qquad |w|<1.
\]
Since $\phi(w)=L w+O(w^2)$ at the origin, the coefficient of $w^{n+1}$ in $G(w)$ is
\[
L^{n+1}\bigl(\tau_{n+1}(\omega')-\tau_{n+1}(\omega)\bigr).
\]
Applying Cauchy's estimate on $|w|=r$ gives
\[
L^{n+1}|\tau_{n+1}(\omega')-\tau_{n+1}(\omega)|
\le 2M r^{-(n+1)}.
\]
Therefore
\[
|\tau_{n+1}(\omega')-\tau_{n+1}(\omega)|
\le \frac{2M}{L r}(L r)^{-n}.
\]
Since $L r>1$, this has the required exponential form with $\vartheta=(L r)^{-1}$.
\end{proof}

\subsection{Resampling future coordinates}

Let $\Omega^+:=\widetilde\Omega^{\mathbb N_0}$, equipped with the product law $\widehat {\mathbb{P}}:=\widetilde {\mathbb{P}}^{\mathbb N_0}$.
The operation $\mathcal T_n$ changes only one symbol.
The operation below changes all future coordinates; nevertheless, $\tau_{n+1}$ only sees the first changed symbol.
For $n\ge0$, $\omega\in\Omega$, and $\eta=(\eta_0,\eta_1,\dots)\in\Omega^+$, define the sequence obtained by resampling the future coordinates from time $n$ by
\[
\widehat{\mathcal R}_n(\omega,\eta)_j=
\begin{cases}
\omega_j,& j<n,\\
\eta_{j-n},& j\ge n.
\end{cases}
\]

\begin{lem}\label{lem:resample-set}
Let $E\subset\Omega$ be measurable.
Then for $\PP$-almost every $\omega$,
\[
\widehat P\bigl\{\eta\in\Omega^+:\widehat{\mathcal R}_n(\omega,\eta)\in E\bigr\}\longrightarrow 1_E(\omega)
\qquad (n\to\infty).
\]
In particular, for $\PP$-almost every $\omega\in E$ one has
\[
\widehat{\mathbb{P}}\bigl\{\eta\in\Omega^+:\widehat{\mathcal R}_n(\omega,\eta)\in E\bigr\}\to1.
\]
\end{lem}

\begin{proof}
Let $\mathcal G_n$ be the $\sigma$-algebra generated by the coordinates $\omega_j$ with $j<n$.
For each bounded measurable function $\Phi$ on $\Omega$,
\[
\int_{\Omega^+}\Phi(\widehat{\mathcal R}_n(\omega,\eta))\,d\widehat{\mathbb{P}}(\eta)
=\EE(\Phi\mid\mathcal G_n)(\omega)
\]
for $\PP$-almost every $\omega$.
Indeed, this identity is immediate for cylinder functions and then follows by the monotone class theorem.
Intersecting the exceptional null sets over the countable set of indices $n$, we may assume that it holds for all $n$ simultaneously.
Applying this with $\Phi=1_E$ gives
\[\widehat{\mathbb{P}}\bigl\{\eta:\widehat{\mathcal R}_n(\omega,\eta)\in E\bigr\}
=\EE(1_E\mid\mathcal G_n)(\omega)
\]
for $\PP$-almost every $\omega$.
The filtration $(\mathcal G_n)_{n\ge0}$ is increasing and generates the full product $\sigma$-algebra, since every fixed coordinate belongs to $\mathcal G_n$ for all sufficiently large $n$.
By the martingale convergence theorem,
\[
\EE(1_E\mid\mathcal G_n)\to 1_E
\qquad\text{almost surely and in }L^1.
\]
This proves the lemma.
\end{proof}

\subsection{The lower bound from one resampling}

For $\omega\in\Omega$, $n\ge1$ and $\eta\in\Omega^+$, define
\[
\Delta_n(\omega,\eta):=\tau_{n+1}(\widehat{\mathcal R}_n(\omega,\eta))-\tau_{n+1}(\omega).
\]

\begin{lem}\label{lem:lower-bound}
Under Assumption \ref{ass:SG2}, for $\PP$-almost every $\omega$, every $n\ge1$ and $\widehat P$-almost every $\eta\in\Omega^+$,
\[
\Delta_n(\omega,\eta)
=\lambda_\omega^n\nu_\omega\bigl((\LL_{\eta_0}-\LL_{\omega_n})h_{\sigma^n\omega}\bigr)+\mathcal E_n(\omega,\eta),
\]
with
\[
|\mathcal E_n(\omega,\eta)|\le C_3\rho^n\lambda_\omega^n.
\]
Consequently, for every $a>0$, for $\PP$-almost every $\omega$ and all sufficiently large $n$ with $\lambda_\omega^n\ge a$,
\[
\widetilde{\mathbb P }\{\eta\in\Omega^+:|\Delta_n(\omega,\eta)|\ge a\delta\}\ge p_0.
\]
\end{lem}

\begin{proof}
The coefficient $\tau_{n+1}$ depends only on the symbols $\omega_0,\dots,\omega_n$.
Therefore
\[
\Delta_n(\omega,\eta)=\tau_{n+1}(\mathcal T_n(\omega,\eta_0))-\tau_{n+1}(\omega), 
\]
where $\mathcal T_n$ was given as (\ref{Tn}).

The formula in Proposition~\ref{prop:sd-consequences}, which uses \hyperref[ass:SGb]{(SGb)}, gives the stated formula and error bound. Since $n$ is countable, the exceptional sets in the future variable may be chosen so that the formula holds for every $n$ on a $\widetilde{\mathbb P}$-full set depending on $\omega$.
Choose $n$ so large that $C_3\rho^n\le\delta$.
Fix $a>0$. If $\lambda_\omega^n\ge a$ and
\[
\bigl|\nu_\omega\bigl((\LL_{\eta_0}-\LL_{\omega_n})h_{\sigma^n\omega}\bigr)\bigr|\ge 3\delta,
\]
then
\[
|\Delta_n(\omega,\eta)|
\ge \lambda_\omega^n(3\delta-C_3\rho^n)
\ge 2a\delta
\ge a\delta.
\]
The last claim follows from item (2) in \hyperref[ass:SG2]{(SG2)},
because the first coordinate $\eta_0$ of the future sample has law $\widetilde{\mathbb{P}}$.
\end{proof}

\subsection{Recurrence of the zero-mean cocycle}

\begin{lem}\label{lem:oscillation}
Let $S$ be an ergodic measure-preserving transformation of the probability space $(X,\mathcal A,\mu)$, and let $f\in L^1(\mu)$ satisfy $\int f\,d\mu=0$.
Then, for every $\varepsilon>0$ and for $\mu$-almost every $x$,
\[
\left|\sum_{j=0}^{n-1} f(S^j x)\right|<\varepsilon
\]
for infinitely many $n$.
\end{lem}

\begin{proof}
This is Atkinson's recurrence theorem for integrable zero-mean cocycles; see Atkinson~\cite{Atk76}.
\end{proof}

\subsection{Proof of the natural-boundary theorem}

\begin{proof}[Proof of Theorem~\ref{thm:boundary}]

Assume, toward contradiction, that $\zeta_\omega$ fails to have the unit circle as a natural boundary on a set of positive probability. Then $\zeta_{\omega}$ extends meromorphically across at least one boundary point on a set of positive probability. By \hyperref[ass:SG2]{(SG2)}, together with that $\EE[\log\lambda]=0$ and Birkhoff's theorem give $\lambda_\omega^n=e^{o(n)}$, $H_\omega$ converge there, which yields $d_{0,\omega}$ extends meromorphically across at least one boundary point on a set of positive probability.
By Lemma~\ref{lem:coeff-decay} there exist a set $\Omega_*$ of positive probability and constants $C>0$ and $\vartheta\in(0,1)$ such that, whenever $\omega,\omega'\in\Omega_*$ satisfy $\tau_j(\omega)=\tau_j(\omega')$ for $1\le j\le n$, we have $|\tau_{n+1}(\omega)-\tau_{n+1}(\omega')|\le C\vartheta^n$.

Let
\[
\Omega_0:=\Omega_*\cap\Omega_1\cap \Omega_2\cap\Omega_3\cap\Omega_4.
\]
Since \(\Omega_0\) has positive probability, we may fix \(\omega\in\Omega_0\) for which the conclusions of Lemma~\ref{lem:resample-set} (with \(E=\Omega_0\)) and Lemmas~\ref{lem:lower-bound} and \ref{lem:oscillation} all hold.
Set
\[
S_n:=\log \lambda_\omega^n=\sum_{j=0}^{n-1}\log\lambda_{\sigma^j\omega}.
\]
By Lemma~\ref{lem:oscillation}, applied to the ergodic Bernoulli shift with $f(\omega)=\log\lambda_\omega$, there are infinitely many $n$ with $\lambda_\omega^n\ge e^{-1}$.
On the auxiliary product probability space $((\Omega^+)^{\mathbb N_0},\overline P)$, where $\overline P:=\widehat {\mathbb{P}}^{\mathbb N_0}$, choose $\eta^{(n)}\in\Omega^+$, each with law $\widehat{\mathbb P}$.
For each such $n$ define
\[
A_n:=\{\widehat{\mathcal R}_n(\omega,\eta^{(n)})\in\Omega_0\},
\qquad
B_n:=\{|\Delta_n(\omega,\eta^{(n)})|\ge e^{-1}\delta\}.
\]
Since $\overline P(A_n)\to1$ by Lemma~\ref{lem:resample-set}, and $\overline P(B_n)\ge p_0$ for all sufficiently large $n$ with $\lambda_\omega^n\ge e^{-1}$ by Lemma~\ref{lem:lower-bound},
\[
\overline P(A_n\cap B_n)
\ge \overline P(B_n)-\overline P(A_n^{c})
\ge p_0/2
\]
for all sufficiently large $n$ with $\lambda_\omega^n\ge e^{-1}$.
The index set $\{n:\lambda_\omega^n\ge e^{-1}\}$ is fixed once $\omega$ is fixed, and the events $A_n\cap B_n$ depend on $\eta^{(n)}$.
Hence these events are independent over this index set and
\[
\sum_{n:\,\lambda_\omega^n\ge e^{-1}}\overline P(A_n\cap B_n)=\infty.
\]
By Borel--Cantelli Lemma, for $\overline P$-almost every sequence, infinitely many such $n$ satisfy both $A_n$ and $B_n$.
Choose one realization of the samples for which this conclusion holds, and then choose such an $n$ with $C\vartheta^n<e^{-1}\delta$.

Because $\widehat{\mathcal R}_n(\omega,\eta^{(n)})$ and $\omega$ have the same coordinates before time $n$,
\[
\tau_j(\widehat{\mathcal R}_n(\omega,\eta^{(n)}))=\tau_j(\omega),
\qquad 1\le j\le n,
\]
by locality.
Since $\Omega_0\subset\Omega_*$, applying Lemma~\ref{lem:coeff-decay} to the pair $\omega,\widehat{\mathcal R}_n(\omega,\eta^{(n)})$ gives
\[
|\Delta_n(\omega,\eta^{(n)})|
=|\tau_{n+1}(\widehat{\mathcal R}_n(\omega,\eta^{(n)}))-\tau_{n+1}(\omega)|
\le C\vartheta^n.
\]
On the other hand, $B_n$ says that $|\Delta_n(\omega,\eta^{(n)})|\ge e^{-1}\delta$.
This is impossible by the choice of $n$.
This contradiction proves the theorem.
\end{proof}

\section{Uniformly expanding maps}\label{sec:example}

In this section, we verify \hyperref[ass:SG]{(SG)} and
item (1) in \hyperref[ass:SG2]{(SG2)} for uniformly expanding
families of $C^{r+1}$ maps with uniformly positive $C^r$ weights
with $\mathcal B=C^1(M,\mathbb R)$.
In particular, we verify \hyperref[ass:SG1]{(SG1)} for the geometric potential.

For simplicity, let $M=\mathbb T^d=\mathbb R^d/\mathbb Z^d$ be the $d$-dimensional torus, with normalized Lebesgue measure $m$.\footnote{The arguments below extend to general compact manifolds by using local coordinates and a partition of unity; cf.~\cite[Sections 2.2 and 3.2.2]{Bal18} and \cite[Section 1.1]{NS19}.}
Let $T$ and $g$ be as in Section~\ref{subsec:basic}.
Without loss of generality, we assume that the conditions on $T_\xi$ in Section~\ref{subsec:basic} hold for every $\xi\in\widetilde\Omega$.
Suppose that, for some $\kappa>1$,
\begin{equation}\label{eq:s5-geometry}
 \sup_{\xi\in\widetilde\Omega}\sup_{x\in M}\|(D_xT_\xi)^{-1}\|\le\kappa^{-1},
 \qquad \inf_{\xi\in\widetilde\Omega}\inf_{x\in M}g_\xi(x)>0.
\end{equation}
We use the norm $\|u\|_{C^1}=\|u\|_{C^0}+\|Du\|_{C^0}$.
For the geometric potential, the lower bound in \eqref{eq:s5-geometry} follows from the uniform bound on $DT_\xi$, and we will obtain \hyperref[ass:SG1]{(SG1)}.
Moreover, Lemma~\ref{lem:s5-spectral} below shows that
$\log\lambda$ is bounded, giving item (3) in \hyperref[ass:SG2]{(SG2)}.
Thus item (2) in \hyperref[ass:SG2]{(SG2)} is the only additional
condition needed to apply Theorem~\ref{thm:boundary} in this setting.
A sufficient condition for it is given in
Proposition~\ref{prop:anti-concentration}.
In what follows, $C$ may change from line to line but is independent of $\omega$,
$n$, $\xi$ and $\varepsilon$.

\subsection{Preliminary estimates and spectral decomposition}\label{subsec:s5-spectral}
The proofs of the following two lemmas are standard (cf. \cite[Sections 2.1--2.3]{Bal18}
and \cite[Section 4]{SR21}).

\begin{lem}\label{lem:s5-bounded}
There is $C>0$ such that
\[
 \sup_{\xi\in\widetilde\Omega}\|\mathcal L_\xi\|_{C^1\to C^1}\le C,\qquad
 \sup_{\xi\in\widetilde\Omega}\|\mathcal L_\xi\|_{C^0\to C^0}\le C.
\]
\end{lem}

\begin{lem}\label{lem:s5-branches}
There is $C>0$ such that
\[
 \|DF\|_{C^0}\le\kappa^{-n},\qquad
 \|D^2F\|_{C^0}\le C\kappa^{-n},\qquad
 \|D\log(g_\omega^{(n)}\circ F)\|_{C^0}\le C,
\]
for every $n\ge1$ and every inverse branch $F$ of $T_\omega^n$.
\end{lem}

We obtain the following spectral decomposition by adapting
the cone argument in \cite[Sections 2 and 4]{SR21}.

\begin{lem}\label{lem:s5-spectral}
There exist measurable families $h_\omega\in C^1(M,\mathbb R)$,
$\nu_\omega\in (C^1(M,\mathbb R))^*$ and $\lambda_\omega>0$,
with $h_\omega$ strictly positive and $\nu_\omega$ positive,
and constants $C>0$ and $\theta\in(0,1)$ such that,
for every $\omega\in\Omega$ and $n\ge1$,
\begin{equation}\label{eq:s5-leading}
 \mathcal L_\omega h_\omega=\lambda_\omega h_{\sigma\omega},
 \quad
 \nu_\omega(h_\omega)=\int_M h_\omega\,dm=1,
 \quad
 \|h_\omega\|_{C^1}+\|\nu_\omega\|_{(C^1)^*}\le C,
\end{equation}
and
\begin{equation}\label{eq:s5-split}
 \mathcal L_\omega^n
 =\lambda_\omega^n h_{\sigma^n\omega}\otimes\nu_\omega+R_\omega^n,
 \qquad
 \|R_\omega^n\|_{C^1\to C^1}\le C\theta^n\lambda_\omega^n.
\end{equation}
Moreover, $h_\omega$ and $\lambda_\omega$ are uniformly
bounded away from zero, and $\lambda_\omega$ is uniformly bounded above.
\end{lem}

\begin{proof}
Write $\ell(u)=\int_M u\,dm$.
For $a>0$, let
\[
 \mathcal C_a=\{u\in C^1:u\ge0,\ |Du(x)|\le au(x)
 \text{ for all }x\in M\}.
\]
By the chain rule, $\mathcal L_\xi(\mathcal C_a)\subset\mathcal C_b$,
where
\[
 b=\kappa^{-1}\bigl(a+\sup_\xi\|D\log g_\xi\|_{C^0}\bigr)
\]
(cf.~\cite[proof of Theorem~4.17]{SR21}
and \cite[Proposition~2.4]{Via97}).
Taking $a$ sufficiently large, we have $b<a$.
Since $\sup u\le e^{a\operatorname{diam}(M)}\inf u$
for $u\in\mathcal C_a$, we have $\|u\|_{C^1}\le C\ell(u)$.
Moreover, since $b<a$, the set
$\{u\in\mathcal C_b:\ell(u)=1\}$ is uniformly bounded away
from the boundary of $\mathcal C_a$ (cf.~\cite[Proposition~2.5]{Via97} and \cite[Lemma~4.16]{SR21}).
Thus $\mathcal C_b\setminus\{0\}$ has finite projective diameter
$\Delta$ in $\mathcal C_a$
(cf.~\cite[Lemma~4.16 and Corollary~A.5]{SR21}
and \cite[Proposition~2.5]{Via97}).

By \cite[Lemma~2.1 and Corollary~A.8]{SR21}, the normalized pullback iterates
$
\ell(\mathcal L_{\sigma^{-n}\omega}^n1)^{-1}\mathcal L_{\sigma^{-n}\omega}^n1
$
converge in $C^1$, uniformly in $\omega$, to a function
$h_\omega\in\mathcal C_b$ with $\ell(h_\omega)=1$ and
$\mathcal L_\omega h_\omega=\lambda_\omega h_{\sigma\omega}$,
where $\lambda_\omega=\ell(\mathcal L_\omega h_\omega)$.
In particular, $h_\omega$ is uniformly bounded in $C^1$
and away from zero. By Lemma~\ref{lem:s5-bounded} and
\eqref{eq:s5-geometry}, $\lambda_\omega$ is uniformly bounded
above and away from zero.

Set $\theta=\tanh(\Delta/4)$ and
$P_\omega^n=(\lambda_\omega^n)^{-1}\mathcal L_\omega^n$.
Differentiating the normalized iterates at $h_\omega$
and using \cite[(2.5)]{SR21}, we obtain
\[
 \|P_\omega^n u-\ell(P_\omega^n u)h_{\sigma^n\omega}\|_{C^1}
 \le C\theta^n\|u\|_{C^1}.
\]
Since $P_\omega h_\omega=h_{\sigma\omega}$ and
$\sup_\omega\|P_\omega\|_{C^1\to C^1}<\infty$, it follows that
\begin{align*}
 |\ell(P_\omega^{n+1}u)-\ell(P_\omega^n u)|
 =|\ell(P_{\sigma^n\omega}
 (P_\omega^n u-\ell(P_\omega^n u)h_{\sigma^n\omega}))| 
 \le C\theta^n\|u\|_{C^1}.
\end{align*}
Thus $(\ell(P_\omega^n u))_{n\ge1}$ is a Cauchy sequence.
Defining $\nu_\omega(u)=\lim_{n\to\infty}\ell(P_\omega^n u)$
and summing the preceding estimate, we obtain
\[
 |\nu_\omega(u)-\ell(P_\omega^n u)|
 \le C\theta^n\|u\|_{C^1}.
\]
This gives \eqref{eq:s5-split}. Moreover, $\nu_\omega$ is
positive, $\nu_\omega(h_\omega)=1$ and
$\sup_\omega\|\nu_\omega\|_{(C^1)^*}<\infty$.

Finally, the derivatives of $T_\xi$ and $g_\xi$ are measurable as pointwise limits of difference quotients.
By separability of the respective $C^1$ spaces, this and the joint measurability of $T_\xi$ and $g_\xi$ imply that $\xi\mapsto T_\xi$ and $\xi\mapsto g_\xi$ are measurable in the $C^1$ topologies
(cf.~\cite[p.~21]{Cra02}).
For each fixed $u\in C^1$, the transfer operator applied to $u$ depends continuously on the map and the weight in the respective $C^1$ topologies, since local inverse branches depend continuously on the map in the $C^1$ topology.
Thus $\omega\mapsto\mathcal L_\omega$ is strongly measurable.
By \cite[Lemma~A.5]{GTQ14} and the definition of $h_\omega$, the map $\omega\mapsto h_\omega$ is measurable in $C^1$.
The measurability of $\lambda_\omega=\ell(\mathcal L_\omega h_\omega)$ then follows from \cite[Lemma~A.6]{GTQ14}.
By change of variables, the functionals $\ell\circ P_\omega^n$ are given by the densities $(\lambda_\omega^n)^{-1}g_\omega^{(n)}|\det DT_\omega^n|$, which are measurable as $L^1$-valued maps.
Since $\ell\circ P_\omega^n$ converges to $\nu_\omega$ in the dual norm, $\nu_\omega$ is also measurable.
\end{proof}

For the geometric potential $g_\xi=|\det DT_\xi|^{-1}$,
we have $\ell(\mathcal L_\xi u)=\ell(u)$ by the change of variables, where $\ell(u)=\int_M u\,dm$.
Thus $\lambda_\omega=1$ and $\nu_\omega=\ell$ by their definitions.
To obtain \hyperref[ass:SG1]{(SG1)}, it therefore suffices
to verify \hyperref[ass:SG]{(SG)}.

By \eqref{eq:s5-leading} and \eqref{eq:s5-split}, we have
\begin{equation}\label{eq:s5-mass}
 \|\mathcal L_\omega^n1\|_{C^0}
 \le\|\mathcal L_\omega^n1\|_{C^1}
 \le C\lambda_\omega^n.
\end{equation}

\begin{remark}
Lemma~\ref{lem:s5-spectral} alone does not yield the flat-trace
estimates in item (1) of \hyperref[ass:SG]{(SG)} and
item (1) of \hyperref[ass:SG2]{(SG2)}, since an estimate of the form
$|\operatorname{tr}(AB)|\le C_A\|B\|_{C^1\to C^1}$
does not hold in general, even for finite-rank operators $B$
and a fixed bounded operator $A$.\footnote{For $A=\operatorname{Id}$,
take $B_N$ to be convolution with the product Fej\'er kernel
of degree $N$. Then $\|B_N\|_{C^1\to C^1}\le1$, whereas
$\operatorname{tr}(B_N)=(N+1)^d$.}
However, in our setting, both $R_\omega^n$ and
$(\mathcal L_\xi-\mathcal L_{\omega_n})R_\omega^n$
can be expressed as linear combinations of transfer operators
associated with uniformly expanding maps and rank-one operators.
This allows us to estimate the errors in the approximation of
their flat traces by regularized flat traces, uniformly in
$\omega$, $n$ and $\xi$ relative to $\lambda_\omega^n$.
Combining these estimates with an operator-norm bound on
the regularized flat trace, we establish the required
flat-trace estimates below.
\end{remark}

\subsection{Regularization of the flat trace}
Choose an even nonnegative function $\chi\in C_c^\infty(\mathbb R^d)$ supported
in the unit ball, with $\int_{\mathbb R^d}\chi(x)\,dx=1$. For
$0<\varepsilon<\varepsilon_0<1/4$, let
\[
 \chi_\varepsilon(x)=\sum_{j\in\mathbb Z^d}\varepsilon^{-d}
 \chi\big((x+j)/\varepsilon\big),\qquad
 J_\varepsilon u(x)=\int_M\chi_\varepsilon(x-y)u(y)\,dm(y).
\]
We take $I_\varepsilon=J_\varepsilon$ in Section~\ref{subsec:basic}. Its kernel
$K_\varepsilon(x,y)=\chi_\varepsilon(x-y)$ is nonnegative, smooth and supported
in $\{(x,y):\operatorname{dist}(x,y)\le\varepsilon\}$.

The following estimates will be used in Section~\ref{subsec:s5-verification}. For every bounded operator
$S:C^1(M,\mathbb R)\to C^1(M,\mathbb R)$,
\begin{equation}\label{eq:s5-regtrace}
 |\operatorname{tr}^{\flat}_\varepsilon(S)|\le C\varepsilon^{-d}\|S\|_{C^1\to C^1}.
\end{equation}
For $h\in C^1(M,\mathbb R)$ and a positive bounded linear functional $\nu$ on $C^1(M,\mathbb R)$, we have
$\operatorname{tr}^{\flat}(h\otimes\nu)=\nu(h)$ and
\begin{equation}\label{eq:s5-rankone}
 |\operatorname{tr}^{\flat}_\varepsilon(h\otimes\nu)-\nu(h)|
 \le C\varepsilon\|h\|_{C^1}\nu(1).
\end{equation}

Since $\chi_\varepsilon$ is nonnegative with integral one and
$DJ_\varepsilon u=J_\varepsilon Du$, we have
$\|J_\varepsilon u\|_{C^1}\le\|u\|_{C^1}$.
By compactness of $M$, both $u$ and $Du$ are uniformly continuous,
so $J_\varepsilon u\to u$ and $J_\varepsilon Du\to Du$ uniformly.
Hence $\|J_\varepsilon u-u\|_{C^1}\to0$ for every
$u\in C^1(M,\mathbb R)$.
Thus all the conditions on $I_\varepsilon$ in
Section~\ref{subsec:basic} are satisfied by taking
$I_\varepsilon=J_\varepsilon$.
Since $J_\varepsilon u\in C^\infty(M,\mathbb R)$,
this also shows that $C^\infty(M,\mathbb R)$ is dense
in $C^1(M,\mathbb R)$.
Condition~\eqref{T} follows from the continuous inclusion
$C^1\hookrightarrow C^0$.

To prove \eqref{eq:s5-regtrace}, we use Fourier series.\footnote{Recall from Section~\ref{subsec:basic}
that $\operatorname{tr}^{\flat}_\varepsilon(S)=\int_M (SK_\varepsilon(\cdot,x))(x)\,dm(x)$
and $\operatorname{tr}^{\flat}(S)=\lim_{\varepsilon\to0}\operatorname{tr}^{\flat}_\varepsilon(S)$
whenever the limit exists.}
Let $e_k(x)=e^{2\pi i k\cdot x}$ for $k\in\mathbb Z^d$ and write
$\widehat u(k)=\int_M u(x)\overline{e_k(x)}\,dm(x)$ and
$\widehat\chi(\zeta)=\int_{\mathbb R^d}\chi(x)e^{-2\pi i\zeta\cdot x}\,dx$.
By the change of variables $z=(x+j)/\varepsilon$ in each term of the sum defining
$\chi_\varepsilon$,
\[
 \widehat{\chi_\varepsilon}(k)
 =\int_{\mathbb R^d}\chi(z)e^{-2\pi i\varepsilon k\cdot z}\,dz
 =\widehat\chi(\varepsilon k).
\]

Fix an integer $N>d/2$.
Since $\chi\in C_c^\infty(\mathbb R^d)$, we have
$|\widehat\chi(\zeta)|\le C_N(1+|\zeta|)^{-2N}$
for some $C_N>0$.
There are at most $C_d\varepsilon^{-d}$ lattice points with $|k|<\varepsilon^{-1}$,
and at most $C_d(2^{j+1}/\varepsilon)^d$ lattice points with
$2^j\varepsilon^{-1}\le |k|<2^{j+1}\varepsilon^{-1}$, $j\ge0$. Hence
\begin{equation}\label{eq:20261128}
 \sum_{k\in\mathbb Z^d}|\widehat\chi(\varepsilon k)|
 \le C\varepsilon^{-d}\left(1+\sum_{j\ge0}2^{j(d-2N)}\right)
 \le C\varepsilon^{-d}.
\end{equation}
By the rapid decay of $\widehat\chi$, the series
\[
 K_\varepsilon(\cdot,y)
 =\sum_{k\in\mathbb Z^d}\widehat\chi(\varepsilon k)e_k(\cdot)\overline{e_k(y)}
\]
converges in $C^1$, uniformly in $y$.
By the boundedness of $S$ and the inclusion $C^1\hookrightarrow C^0$, we may interchange the sum
and the integral to obtain
\[
 \operatorname{tr}^{\flat}_\varepsilon(S)
 =\sum_{k\in\mathbb Z^d}\widehat\chi(\varepsilon k)
       \int_M(Se_k)(x)\overline{e_k(x)}\,dm(x)
 =\sum_{k\in\mathbb Z^d}\widehat\chi(\varepsilon k)\widehat{Se_k}(k).
\]
For $k\ne0$, integrating by parts in the direction $k/|k|$ gives
\[
 |\widehat{Se_k}(k)|\le\frac{\|D(Se_k)\|_{C^0}}{2\pi|k|}
 \le C\|S\|_{C^1\to C^1},
\]
since $\|e_k\|_{C^1}=1+2\pi|k|$. The same bound holds for $k=0$ by the definition of the Fourier coefficient.
Thus $|\operatorname{tr}^{\flat}_\varepsilon(S)|\le C\|S\|_{C^1\to C^1}\sum_k|\widehat\chi(\varepsilon k)|$,
which gives \eqref{eq:s5-regtrace} by \eqref{eq:20261128}.

For \eqref{eq:s5-rankone}, we first note that
\[
 \operatorname{tr}^{\flat}_\varepsilon(h\otimes\nu)
 =\int_M h(x)\nu(K_\varepsilon(\cdot,x))\,dm(x)=\nu(J_\varepsilon h),
\]
where the last equality follows from the continuity of $\nu$ on $C^1$.
In particular, $\operatorname{tr}^{\flat}(h\otimes\nu)=\nu(h)$.
By positivity of $\nu$ and the mean value theorem,
\[
 |\nu(J_\varepsilon h-h)|\le\nu(1)\|J_\varepsilon h-h\|_{C^0}
 \le\varepsilon\nu(1)\|Dh\|_{C^0}\int_{\mathbb R^d}|z|\chi(z)\,dz,
\]
which gives \eqref{eq:s5-rankone}.

\subsection{Verification of \texorpdfstring{\hyperref[ass:SG]{(SG)}}{(SG)} and item (1) in \texorpdfstring{\hyperref[ass:SG2]{(SG2)}}{(SG2)}}
\label{subsec:s5-verification}

We first establish a uniform rate of convergence for the
regularized flat trace of $\mathcal L_\omega^n$.

\begin{lem}\label{lem:s5-quantitative}
There is $C>0$ such that, for every $\omega\in\Omega$ and $n\ge1$,
the flat trace of $\mathcal L_\omega^n$ exists and
\begin{equation}\label{eq:20261128b}
 \operatorname{tr}^{\flat}(\mathcal L_\omega^n)
 =\sum_{T_\omega^n y=y}
 \frac{g_\omega^{(n)}(y)}
 {|\det(\operatorname{Id}-(D_yT_\omega^n)^{-1})|}.
\end{equation}
Moreover, for every $\varepsilon\in(0,\varepsilon_0)$,
\begin{equation}\label{eq:s5-quantitative}
 \left|\operatorname{tr}^{\flat}_\varepsilon(\mathcal L_\omega^n)
 -\operatorname{tr}^{\flat}(\mathcal L_\omega^n)\right|
 \le C\varepsilon\|\mathcal L_\omega^n1\|_{C^0}.
\end{equation}
\end{lem}

\begin{proof}
We follow the 
argument in the proof of
\cite[Proposition 3.14]{Bal18}.
Fix $\omega\in\Omega$ and $n\ge1$, and define
$\Psi:M\to M$ by $\Psi(y)=T_\omega^n(y)-y$ modulo $\mathbb Z^d$.
Since $\|(D_yT_\omega^n)^{-1}\|\le\kappa^{-n}<1$,
$D_y\Psi$ is invertible for every $y\in M$.
By compactness and connectedness, $\Psi$ is a finite covering
of $M$. We can thus define the positive $C^1$ function
\[
 b(w)=\sum_{\Psi(y)=w}
 \frac{g_\omega^{(n)}(y)}
 {|\det(\operatorname{Id}-(D_yT_\omega^n)^{-1})|}.
\]
By changing variables first for $T_\omega^n$
and then for $\Psi$, we have
\begin{align*}
 \operatorname{tr}^{\flat}_\varepsilon(\mathcal L_\omega^n)
 &=\int_M\sum_{T_\omega^n y=x}
 g_\omega^{(n)}(y)\chi_\varepsilon(y-x)\,dm(x)\\
 &=\int_M g_\omega^{(n)}(y)|\det D_yT_\omega^n|
 \chi_\varepsilon(-\Psi(y))\,dm(y)\\
 &=\int_M\chi_\varepsilon(-w)b(w)\,dm(w),
\end{align*}
where we used
$
 |\det D_yT_\omega^n|/|\det D_y\Psi|
 =1/|\det(\operatorname{Id}-(D_yT_\omega^n)^{-1})|.
$
By the definition of $\chi_\varepsilon$ and the mean value theorem,
\begin{align*}
 \left|\operatorname{tr}^{\flat}_\varepsilon(\mathcal L_\omega^n) -b(0)\right|
 =\left|\int_{\mathbb R^d}\chi(z)(b(-\varepsilon z)-b(0))\,dz\right| 
 \le\varepsilon\|Db\|_{C^0} \int_{\mathbb R^d}|z|\chi(z)\,dz,
\end{align*}
where $b$ is viewed as a $\mathbb Z^d$-periodic function on $\mathbb R^d$.
Letting $\varepsilon\to0$ gives \eqref{eq:20261128b}, since
$\Psi(y)=0$ if and only if $T_\omega^n(y)=y$.
Thus, to prove \eqref{eq:s5-quantitative}, it remains to show that
\[
 \|Db\|_{C^0}\le C\|\mathcal L_\omega^n1\|_{C^0}.
\]

For each $y\in M$, let $F$ be an inverse branch of $T_\omega^n$
defined near $T_\omega^n(y)$ with $F(T_\omega^n(y))=y$.
Since $\|DF\|_{C^0}\le\kappa^{-n}$,
the map $\operatorname{Id}-F$ is locally invertible and,
by the triangle inequality,
\[
 \|D((\operatorname{Id}-F)^{-1})\|_{C^0}
 \le (1-\|DF\|_{C^0})^{-1}
 \le (1-\kappa^{-1})^{-1}.
\]
The identity $\Psi\circ F=\operatorname{Id}-F$
implies that the corresponding inverse branch of $\Psi$
is given by $F\circ(\operatorname{Id}-F)^{-1}$.
Thus, the corresponding summand in $b$ is locally given by
\[
 \left(
 \frac{g_\omega^{(n)}\circ F}
 {|\det(\operatorname{Id}-DF)|}
 \right)\circ(\operatorname{Id}-F)^{-1}.
\]
For each unit vector $v$, we have
\[
 \partial_v\log|\det(\operatorname{Id}-DF)|
 =-\operatorname{tr}\bigl(
 (\operatorname{Id}-DF)^{-1}\partial_v DF\bigr).
\]
By Lemma~\ref{lem:s5-branches} and the chain rule,
the norm of the logarithmic derivative of each summand in $b$
is therefore bounded by
$
 \bigl(C+dC\kappa^{-1}(1-\kappa^{-1})^{-1}\bigr)
 (1-\kappa^{-1})^{-1}.
$
By positivity, $\|D_wb\|\le Cb(w)$ for all $w\in M$.

For any $w_1,w_2\in M$, let $\gamma:[0,\ell]\to M$
be a shortest geodesic from $w_1$ to $w_2$ parametrized
by arc length. Since $\|D\log b\|_{C^0}\le C$, we have
\begin{align*}
 |\log b(w_2)-\log b(w_1)|
 =\left|\int_0^\ell
 D_{\gamma(t)}\log b\,(\dot\gamma(t))\,dt\right|
 \le C\ell
 \le C\operatorname{diam}(M).
\end{align*}
Consequently, $\sup b\le e^{C\operatorname{diam}(M)}\inf b$.
On the other hand, by changing variables we obtain
\[
 \int_M b(w)\,dm(w)
 =\int_M g_\omega^{(n)}(y)|\det D_yT_\omega^n|\,dm(y)
 =\int_M\mathcal L_\omega^n1\,dm.
\]
Since $m(M)=1$, we obtain
\[
 \|Db\|_{C^0}\le C\sup b
 \le C\int_M b\,dm
 \le C\|\mathcal L_\omega^n1\|_{C^0},
\]
as desired.
\end{proof}

We now verify item (1) in \hyperref[ass:SG]{(SG)} and
item (1) in \hyperref[ass:SG2]{(SG2)}.
By Lemma~\ref{lem:s5-bounded}, \eqref{eq:s5-leading}
and \eqref{eq:s5-split},
\[
 \|(\mathcal L_\xi-\mathcal L_{\omega_n})R_\omega^n\|_{C^1\to C^1}
 \le C\theta^n\lambda_\omega^n,
 \qquad
 \|(\mathcal L_\xi-\mathcal L_{\omega_n})
 h_{\sigma^n\omega}\|_{C^1}\le C.
\]
Applying \eqref{eq:s5-quantitative} to the two compositions in
\[
 (\mathcal L_\xi-\mathcal L_{\omega_n})\mathcal L_\omega^n
 =\mathcal L_\xi\mathcal L_\omega^n-\mathcal L_\omega^{n+1},
\]
we obtain
\begin{align*}
 &\left|
 \operatorname{tr}^{\flat}
 ((\mathcal L_\xi-\mathcal L_{\omega_n})\mathcal L_\omega^n)
 -\operatorname{tr}^{\flat}_\varepsilon
 ((\mathcal L_\xi-\mathcal L_{\omega_n})\mathcal L_\omega^n)
 \right| \le C\varepsilon\lambda_\omega^n,
\end{align*}
since
$
 \|\mathcal L_\xi\mathcal L_\omega^n1\|_{C^0}
 +\|\mathcal L_\omega^{n+1}1\|_{C^0}
 \le C\lambda_\omega^n
$
by Lemma~\ref{lem:s5-bounded} and \eqref{eq:s5-mass}.
Since
\[
 (\mathcal L_\xi-\mathcal L_{\omega_n})R_\omega^n
 =(\mathcal L_\xi-\mathcal L_{\omega_n})\mathcal L_\omega^n
 -\lambda_\omega^n
 ((\mathcal L_\xi-\mathcal L_{\omega_n})h_{\sigma^n\omega})
 \otimes\nu_\omega,
\]
its flat trace exists.
By the preceding estimates, the uniform bound on $\nu_\omega(1)$ from \eqref{eq:s5-leading}, and \eqref{eq:s5-rankone},
\begin{align*}
 &\left|
 \operatorname{tr}^{\flat}
 ((\mathcal L_\xi-\mathcal L_{\omega_n})R_\omega^n)
 -\operatorname{tr}^{\flat}_\varepsilon
 ((\mathcal L_\xi-\mathcal L_{\omega_n})R_\omega^n)
 \right|\le C\varepsilon\lambda_\omega^n.
\end{align*}
Combining this with \eqref{eq:s5-regtrace}, we obtain
\[
 |\operatorname{tr}^{\flat}
 ((\mathcal L_\xi-\mathcal L_{\omega_n})R_\omega^n)|
 \le C(\varepsilon+\varepsilon^{-d}\theta^n)\lambda_\omega^n.
\]
Taking $\rho=\theta^{1/(d+1)}$ and $\varepsilon=\varepsilon_0\rho^n$
yields
\[
 |\operatorname{tr}^{\flat}
 ((\mathcal L_\xi-\mathcal L_{\omega_n})R_\omega^n)|
 \le C\rho^n\lambda_\omega^n.
\]
The same argument applies to $R_\omega^n$ itself, using
\eqref{eq:s5-leading}, \eqref{eq:s5-split} and
\eqref{eq:s5-quantitative} directly, and gives
\begin{equation}\label{eq:s5-rem}
 |\operatorname{tr}^{\flat}(R_\omega^n)|
 \le C\rho^n\lambda_\omega^n.
\end{equation}
Since $\theta<\rho<1$ and $|\nu_\omega(h_{\sigma^n\omega})|$
is uniformly bounded by \eqref{eq:s5-leading},
this verifies item (1) in \hyperref[ass:SG]{(SG)}
and item (1) in \hyperref[ass:SG2]{(SG2)} with the same $\rho$.

It remains to compare $Z_n(\omega)$ and $\tau_n(\omega)$.
By \eqref{eq:s5-geometry} and \eqref{eq:20261128b},
\[
 |Z_n(\omega)-\tau_n(\omega)|
 \le \tau_n(\omega)\sup_{x\in M}
 \big||\det(\operatorname{Id}-(D_xT_\omega^n)^{-1})|-1\big|
 \le C\kappa^{-n}\tau_n(\omega).
\]
Here we used the inequality
$\big||\det(\operatorname{Id}-A)|-1\big|
\le (2^d-1)\|A\|$ for $\|A\|\le1$,
which follows from the multilinearity of the determinant.
On the other hand, by the identity $\operatorname{tr}^{\flat}(h\otimes\nu)=\nu(h)$, \eqref{eq:s5-leading}
and \eqref{eq:s5-rem}, we have
$\tau_n(\omega)\le C\lambda_\omega^n$.
This gives item (2) in \hyperref[ass:SG]{(SG)} with $\alpha=\kappa$,
completing the verification of \hyperref[ass:SG]{(SG)} and
item (1) in \hyperref[ass:SG2]{(SG2)}.

\subsection*{Acknowledgements}
Y.~Nakano was supported by JSPS KAKENHI Grant Number JP23K03188
and by JST PRESTO Grant Number JPMJPR25K8.
S.~Suzuki was supported by JSPS KAKENHI Grant Number JP24K16932.

\subsection*{AI statement}
The main ideas and overall framework of this paper were developed by the authors.
During the preparation of the manuscript,
ChatGPT was used primarily to assist in refining the mathematical exposition, checking mathematical arguments, and language editing. The authors have reviewed the AI-assisted material and take full responsibility for the contents of this paper.

\end{document}